\documentclass[11pt]{article}
\usepackage[english]{babel}
\usepackage{amsmath}
\usepackage{amsthm}
\usepackage{amsfonts}
\usepackage{amssymb}
\usepackage{enumerate}
\usepackage{color}
\usepackage{xspace}
\usepackage{euscript}
\usepackage{graphicx}
\usepackage{amscd}
\usepackage{epsfig}
\usepackage{tabularx}
\usepackage{url}
\usepackage{multimedia}

        \newtheorem{lemma}{Lemma}[section]

         \newtheorem{remark}[lemma]{Remark}
        \newtheorem{definition}{Definition}[section]

\title{A Multilevel Spectral Indicator Method for Eigenvalues of Large Non-Hermitian Matrices}
\author{Ruihao Huang \thanks{Michigan Technological University, Houghton, MI 49931, U.S.A. ({\tt ruihaoh@mtu.edu}).}
\and Jiguang Sun \thanks{Michigan Technological University, Houghton, MI 49931, U.S.A. ({\tt jiguangs@mtu.edu}).}
\and Chao Yang \thanks{ Lawrence Berkeley National Laboratory, Berkeley, CA 94720, U.S.A. ({\tt cyang@lbl.gov}).}}

\begin{document}
\maketitle

\begin{abstract}
Recently a novel family of eigensolvers, called spectral indicator methods (SIMs), was proposed.
Given a region on the complex plane, SIMs first compute an indicator by the spectral projection.
The indicator is used to test if the region contains eigenvalue(s).
Then the region containing eigenvalues(s) is subdivided and tested. The procedure is repeated until the eigenvalues are identified within a specified precision. 
In this paper, using Cayley transformation and Krylov subspaces, 
a memory efficient multilevel eigensolver is proposed. 
The method uses less memory compared with the early versions of SIMs and is particularly suitable to compute many eigenvalues of large sparse (non-Hermitian) matrices. 
Several examples are presented for demonstration.
\end{abstract}

%
%
%
%

\section{Introduction}
Consider the generalized eigenvalue problem
\begin{equation} \label{AxLambdaBx}
A x= \lambda B x,
\end{equation}
where $A, B$ are $n \times n$ large sparse non-Hermitian matrices.
In particular, we are interested in the computation of all
eigenvalues in a region $R \subset \mathbb C$, which contains $p$ eigenvalues such that $1 \ll p \ll n$ or $1 \ll p \sim n$.

Many efficient eigensolvers are proposed in literature for large sparse Hermitian (or symmetric) matrices (see, e.g., \cite{Saad2003}). 
In contrast, for non-Hermitian matrices, there exist much fewer methods including the Arnoldi method
and Jacobi-Davidson method \cite{arpack, Templates2000}. Unfortunately, these methods are still far from
satisfactory as pointed out in \cite{Saad2011}: ``{\it In essence what differentiates the Hermitian from the non-Hermitian eigenvalue problem is
that in the first case we can always manage to compute an approximation whereas there are non-symmetric problems that can be arbitrarily
difficult to solve and can essentially make any algorithm fail}." 

Recently, a family of eigensolvers, called the spectral indicator methods (SIMs), was proposed \cite{SunZhou2016, Huang2016JCP,Huang2018NLA}.
The idea of SIMs is different from the classical eigensolvers. In brief, 
given a region $R \subset \mathbb C$ whose boundary $\partial R$ is a simple closed curve, 
an indicator $I_R$ is defined and then used to decide if $R$ contains eigenvalue(s).
When the answer is positive, $R$ is divided into sub-regions
and indicators for these sub-regions are computed. 
The procedure continues until the size of the sub-region(s) is smaller than the specified precision, e.g., $10^{-6}$.
The indicator $I_R$ is defined using the spectral projection $P$, i.e., Cauchy contour integral of the resolvent 
of the matrix pencil $(A, B)$ on $\partial R$ \cite{Kato1966}.
In particular, one can construct $I_R$ based on the spectral projection of a random vector ${\boldsymbol f}$. 
It is well-known that $P$ projects ${\boldsymbol f}$
to the generalized eigenspace associated to the eigenvalues enclosed by $\partial R$ \cite{Kato1966}. $P{\boldsymbol f}$ is zero if there is no eigenvalue(s)
inside $R$, and nonzero otherwise. Hence $P{\boldsymbol f}$ can be used to decide if $R$ contains eigenvalues(s) or not.
Evaluation of $P{\boldsymbol f}$ needs to solve linear systems at quadrature points on $\partial R$.
In general, it is believed that computing eigenvalues is more difficult than solving linear systems of equations \cite{HernandezEtal2005ACMTOMS}. 
The proposed method converts the eigenvalue problem to solving a number of related linear systems. 

Spectral projection is a classical tool in functional analysis to study, e.g., the spectrum of operators \cite{Kato1966}
and the finite element convergence theory for eigenvalue problems of partial differential equations \cite{SunZhou2016}. 
It has been used to compute matrix eigenvalue problems in the method by Sakurai-Sugiura \cite{SakuraiSugiura2003CAM}
and FEAST by Polizzi \cite{Polizzi2009PRB}. For example, FEAST uses spectral projection to build subspaces and 
thus can be viewed as a subspace method \cite{TangPolizzi2014SIAMMAA}. 
In contrast, SIMs uses the spectral projection to define indicators 
and combines the idea of bisection to locate eigenvalues.
Note that the use of other tools such as the condition number to define the indicator is possible.

In this paper, we propose a new SIM, called SIM-M.
Firstly, by proposing a new indicator, 
the memory requirement is significantly reduced and thus the computation of many eigenvalues of large matrices becomes realistic. 
Secondly, a new strategy to speedup the computation of the indicators is presented. Thirdly, 
other than the recursive calls in the first two members of SIMs \cite{Huang2016JCP,Huang2018NLA}, a multilevel technique
is used to further improve the efficiency. Moreover, a subroutine is added to find the multiplicities of the eigenvalues.
The rest of the paper is organized as follows. Section \ref{CT} presents the basic idea of SIMs and two early members of SIMs.
In Section \ref{MEI}, we propose a new eigensolver SIM-M with the above features. 
The algorithm and the implementation details are discussed as well.
The proposed method is tested by various matrices in Section \ref{NE}. Finally, in Section~\ref{Con}, we draw some conclusions and discuss some future work. 

\section{Spectral Indicator Methods}\label{CT}
In this section, we give an introduction of SIMs and refer the readers to \cite{SunZhou2016, Huang2016JCP,Huang2018NLA} for more details.
For simplicity, assume that $R$ is a square and $\Gamma:=\partial R$ lies in the resolvent set $\mathcal{R}$ of $(A, B)$, i.e., the set of
$z \in \mathbb C$ such that $(A - zB)$ is invertible. The key idea of SIMs is to find an indicator that can be used
to decide if $R$ contains eigenvalue(s) . 

One way to define the indicator is to use the spectral projection, 
a classical tool in functional analysis \cite{Kato1966}.
Specifically, the matrix $P$ defined by
\begin{equation}\label{P}
P=\dfrac{1}{2\pi i}\int_{\Gamma}(A-zB)^{-1}dz
\end{equation}
is the spectral projection of a vector ${\boldsymbol f}$ onto the algebraic eigenspace associated with the eigenvalues of \eqref{AxLambdaBx}
inside $\Gamma$. If there are no eigenvalues inside $\Gamma$, then $P = 0$, and hence $P{\boldsymbol f} = {\boldsymbol 0}$ for all
${\boldsymbol f} \in \mathbb C^n$. If $\Gamma$ does enclose one or more eigenvalues, then $P{\boldsymbol f} \ne {\boldsymbol 0}$ with 
probability $1$ for a random vector ${\boldsymbol f}$.

To improve robustness, in RIM (recursive integral method) \cite{Huang2016JCP}, the first member of SIMs, the indicator is defined as
\begin{equation}\label{indicatorRIM}
I_R:=\left \| P \left( \frac{P {\boldsymbol f}}{\|P {\boldsymbol f}\|}\right)\right\|.
\end{equation}
Analytically, $I_R = 1$ if there exists at least one eigenvalue in $\Gamma$. Note that when a quadrature rule is applied, $I_R \ne 1$ in general. 
The RIM algorithm  is very simple and listed as follows \cite{Huang2016JCP}.
\begin{itemize}
\item[] {\bf RIM}$(A, B, R, h_0, \delta_0, {\boldsymbol f})$
\item[]{\bf Input:}  matrices $A, B$, region $R$, precision $h_0$, threshhold $\delta_0$, random vector ${\boldsymbol f}$.
\item[]{\bf Output:}  generalized eigenvalue(s) $\lambda$ inside $R$
\item[1.] Compute ${I_R}$.
\item[2.] If $I_R < \delta_0$, exit (no eigenvalues in $R$).
\item[3.] Otherwise, compute the diameter $h$ of $R$.
			\begin{itemize}
				\item[-] If $h  > h_0 $, 
				partition $R$ into subregions $R_j, j=1, \ldots N$.
						\begin{itemize}
						\item[] for $j=1$ to $N$
						\item[] $\qquad${\bf RIM}$(A, B, R_j, h_0, \delta_0, {\boldsymbol f})$.
						\item[] end
						\end{itemize}
				\item[-] else, 
						\begin{itemize}
							\item[] set $\lambda$ to be the center of $R$. output $\lambda$ and exit.
						\end{itemize}
			\end{itemize}
	
\end{itemize}

The major task of RIM is to compute the indicator $I_R$ defined in \eqref{indicatorRIM}. Let the approximation to $P{\boldsymbol f}$ be given by
\begin{equation}\label{XLXf}
P{\boldsymbol f} \approx  \dfrac{1}{2 \pi i} \sum_{j=1}^{n_0} \omega_j {\boldsymbol x}_j,
\end{equation}
where $\omega_j$'s are quadrature weights and ${\boldsymbol x}_j$'s are the solutions of the linear systems
\begin{equation}\label{linearsys}
(A- z_jB){\boldsymbol x}_j = {\boldsymbol f}, \quad j = 1, \ldots, n_0.
\end{equation}
Here $z_j's$ are quadrature points on $\Gamma$. 
The total number of the linear systems \eqref{linearsys} for RIM to solve  is at most
\begin{equation}\label{complexityRIM}
 2n_0 \lceil \log_2(h/h_0)\rceil p,
\end{equation}
where $p$ is the number of eigenvalues in $R$, $n_0$ is the number of the quadrature points, $h$ is the size of the $R$, $h_0$ is
the required precision, and $\lceil \cdot \rceil$ denotes the least larger integer. Given $R$, $p$ is a fixed number. 
The complexity of RIM is proportional to the complexity of solving the linear system \eqref{linearsys}.

The computational cost of RIM mainly comes from solving the linear systems \eqref{linearsys} 
to approximate the spectral projection $P{\boldsymbol f}$. 
It is clear that the cost will be greatly reduced if one can take advantage of the parametrized 
linear systems of the same structure. In \cite{Huang2018NLA}, a new member RIM-C (recursive integral method using Cayley transformation) is proposed.
The idea is to construct some Krylov subspaces
and use them to solve \eqref{linearsys} for all quadrature points $z_j$'s.
Since the method we shall propose is based on RIM-C, a description of RIM-C is included as follows.

Let $M$ be a $n \times n$ matrix, ${\boldsymbol b} \in \mathbb C^n$ be a vector, and $m$ be a non-negative integer. 
The Krylov subspace is defined as
\begin{align} 
K_{m}(M; {\boldsymbol b}):=\text{span} \{{\boldsymbol b}, M {\boldsymbol b}, \ldots, M^{m-1}{\boldsymbol b} \}.
\end{align}
It has the shift-invariant property
\begin{align} \label{eq:5}
K_{m}(\gamma_1 M+ \gamma_2 I; {\boldsymbol b})=K_{m}(M; {\boldsymbol b}),
\end{align}
where $\gamma_1$ and $\gamma_2$ are two scalars.

Consider a family of linear systems
\begin{equation} \label{AzBxf}
(A-zB) {\boldsymbol x} ={\boldsymbol f},
\end{equation}
where $z$ is a complex number. 
Assume that $\sigma$ is not a generalized eigenvalue and $\sigma \neq z$. 
By Cayley transformation, multiplying both sides of \eqref{AzBxf} by $(A- \sigma B)^{-1}$, we have that
\begin{eqnarray*}
\label{ABsigmaz}(A- \sigma B)^{-1}{\boldsymbol f}&=&(A- \sigma B)^{-1}(A-z B){\boldsymbol x} \\
\nonumber &=&(A- \sigma B)^{-1}(A-\sigma B+(\sigma -z)B) {\boldsymbol x} \label{eq:Cayley} \\
\nonumber &=&(I+(\sigma - z)(A- \sigma B)^{-1}B) {\boldsymbol x}.
\end{eqnarray*}
Let $M=(A- \sigma B)^{-1}B$ and 
${\boldsymbol b}=(A- \sigma B)^{-1}{\boldsymbol f}$. Then
\eqref{AzBxf} becomes
\begin{align} \label{IMxb}
(I+(\sigma -z)M) {\boldsymbol x} = {\boldsymbol b}.
\end{align}
From \eqref{eq:5}, the Krylov subspace $K_{m}(I+(\sigma -z)M; {\boldsymbol b})$ is the same as $K_{m}(M; {\boldsymbol b})$. 
We shall use $K^\sigma_{m}(M; {\boldsymbol b})$ when it is necessary to indicate its dependence on the shift $\sigma$.

Arnoldi's method is used by RIM-C to solve the linear systems. 
First, consider the orthogonal projection method for 
\[
M {\boldsymbol x}={\boldsymbol b}.
\]
Let the initial guess be ${\boldsymbol x}_0={\boldsymbol 0}$. One seeks an approximate solution ${\boldsymbol x}_m$ in 
$K_m(M; {\boldsymbol b})$ by imposing the Galerkin condition \cite{Saad2003}
\begin{equation} \label{eq:Krylov}
({\boldsymbol b}-M {\boldsymbol x}_m) \perp K_m(M; {\boldsymbol b}).
\end{equation}
The Arnoldi's method (Algorithm 6.1 of \cite{Saad2011}) is as follows.
\begin{itemize}
\item[1.] Choose a vector ${\boldsymbol v}_1$ of norm $1$ (${\boldsymbol v}_1={\boldsymbol b}/\|{\boldsymbol b}\|_2$).
\item[2.] for $j=1,2, \ldots, m$
	\begin{itemize}
		\item $h_{ij} = (M{\boldsymbol v}_j, {\boldsymbol v}_i), \quad i=1,2, \ldots,j$.
		\item ${\boldsymbol w}_j = M {\boldsymbol v}_j - \sum_{i=1}^j h_{ij} {\boldsymbol v}_i$.
		\item $h_{j+1,j}=\|{\boldsymbol w}_j\|_2$. If $h_{j+1, j}=0$, stop.
		\item ${\boldsymbol v}_{j+1} = {\boldsymbol w}_j/h_{j+1, j}$.
	\end{itemize}
\end{itemize}

Let $V_m$ be the $n \times m$ orthogonal matrix with column vectors ${\boldsymbol v}_1, \ldots, {\boldsymbol v}_m$ and $H_m$ be
the $m \times m$ Hessenberg matrix whose nonzero entries are $h_{i,j}$.
Proposition 6.5 of \cite{Saad2011} implies that
\begin{equation} \label{eq:arnoldi}
M V_m=V_m H_m + {\boldsymbol v}_{m+1} {h}_{m+1,m}{\boldsymbol e}^{T}_m
\end{equation}
and
\[
\text{span}\{\mathrm{col}(V_m)\}=K_{m}(M; {\boldsymbol b}).
\] 
Let ${\boldsymbol x}_m=V_m \tilde{\boldsymbol y}$ such that
the Galerkin condition \eqref{eq:Krylov} holds, i.e., 
\begin{equation}
V^{T}_m {\boldsymbol b}-V^{T}_m M V_m \tilde{\boldsymbol y} ={\boldsymbol 0}.
\end{equation}
Using \eqref{eq:arnoldi}, the residual is given by
\begin{equation} \label{eq:err1}
\|{\boldsymbol b}-M {\boldsymbol x}_m\|_2 = {h}_{m+1,m}|{\boldsymbol e}_m^{T} \tilde{\boldsymbol y}|.
\end{equation}

Next, we consider the linear system \eqref{IMxb}. For $I+(\sigma -z)M$, due to the shift invariant property, one has that
\begin{equation} \label{shiftedLS}
\{I+(\sigma -z)M \} V_m =V_m( I+(\sigma -z)H_m)\\
+(\sigma - z){\boldsymbol v}_{m+1} {h}_{m+1,m}{\boldsymbol e}^{T}_m.
\end{equation}
The Galerkin condition \eqref{eq:Krylov} becomes
\begin{equation}
V^{T}_m {\boldsymbol b}-V^{T}_m \{I+(\sigma -z)M \} V_m {\boldsymbol y} =0.
\end{equation}
It implies that
\begin{equation} \label{eq:many}
\{I+(\sigma -z)H_m \} {\boldsymbol y} =\beta {\boldsymbol e}_1,
\end{equation}
where $\beta = \|{\boldsymbol b}\|_2$.
Combination of \eqref{shiftedLS} and \eqref{eq:many} gives the residual
\begin{equation}\label{eq:error}
\|b-\{I+(\sigma - z)M \} {\boldsymbol x}_m\|_2 =(\sigma -z){h}_{m+1,m}|{\boldsymbol e}_m^{T} {\boldsymbol y}|.
\end{equation}
Let $z_j$ be a quadrature point and one need to solve
\begin{align}
(I+(\sigma -z_j)M) {\boldsymbol x}_j={\boldsymbol b},
\end{align}
where $M=(A-\sigma B)^{-1} B$  and  ${\boldsymbol b}=(A- \sigma B )^{-1} {\boldsymbol f}$. 

From \eqref{XLXf} and \eqref{eq:many}, 
\begin{eqnarray} \label{eq:y}
{\boldsymbol y}_j&=&\beta (I+ (\sigma-z_j) H_m)^{-1}{\boldsymbol e}_1,  \\
\nonumber {\boldsymbol x}_j &\approx& V_m {\boldsymbol y}_j, \\
 \label{eq:reduced}
   P{\boldsymbol f} &\approx&  \dfrac{1}{2 \pi i} \sum w_j V_m {\boldsymbol y}_j.
\end{eqnarray}
The idea of RIM-C is to use the Krylov subspace for $M=(A-\sigma B)^{-1} B$ to solve \eqref{linearsys}
for as many $z_j$'s as possible. The residual can be monitored with a little extra cost using \eqref{eq:error}.

Since the Krylov subspace method is used, the indicator defined in \eqref{indicatorRIM} is not appropriate
since it projects ${\boldsymbol f}$ twice.
RIM-C defines an indicator different from \eqref{indicatorRIM}. 
Let $P{\boldsymbol  f}|_{n_0}$ be the approximation of $ P {\boldsymbol f} $ with $n_0$ quadrature points for the circle circumscribing $R$. 
It is well-known that the trapezoidal quadrature
of a periodic function converges exponentially \cite[Section 4.6.5]{davis1984methods}, i.e., 
\begin{align*}
 \left \|P{\boldsymbol f}- P{\boldsymbol f}|_{n_0}\right\| = O(e^{-C n_0}),
  \end{align*} 
where C is a constant. 
For a large enough $n_0$, one has that
\begin{equation}\label{IndicatorPf}
  \dfrac{ \left \| P {\boldsymbol f}|_{2n_0}\right\|}{ \left \| P {\boldsymbol f}|_{n_0} \right\|}=\begin{cases}
  \dfrac{\|P {\boldsymbol f}\| + O(e^{-C 2n_0})}{\|P {\boldsymbol f}\| + O(e^{-C n_0})}  & \text{if there are eigenvalues inside } R, \\
 \dfrac{ O(e^{-C 2n_0})}{ O(e^{-C n_0})}=O(e^{-C n_0})  & \text{no eigenvalue inside } R.
  \end{cases}  
\end{equation}
The indicator is then defined as 
\begin{equation}\label{ISPf}
I_R := \frac{\|P {\boldsymbol f}_{2n_0}\|}{\|P {\boldsymbol f}_{n_0}\|}  
\approx \frac{\quad \left \|\sum_{j=1}^{2n_0} w_j V_m {\boldsymbol y}_j \right\| \quad}{\left\| \sum_{j=1}^{n_0} w_j V_m {\boldsymbol y}_j\right\|}.
\end{equation}

\section{Multilevel Memory Efficient Method}\label{MEI}
In this section, we make several improvements of RIM-C and propose a multilevel memory efficient method, called SIM-M.

\subsection{A New Memory Efficient Indicator}
In view of \eqref{ISPf}, the computation of the indicator needs to store $V_m$. When $R$ contains
a lot of eigenvalues, the method can become memory intensive. 

\begin{definition}
A (square) region $R$ is resolvable with respect to $(\sigma, \epsilon_0)$ if the linear systems \eqref{linearsys} associated with all the quadrature points can be
solved up to the given residual $\epsilon_0$ using the Krylov subspace related to a shift $\sigma$. 
\end{definition}

Assume that $R$ is resolvable with respect to $(\sigma, \epsilon_0)$. From \eqref{ISPf}, one has that
\begin{equation}
    I_R  \approx \frac{ \| \sum_{j=1}^{2n_0} w_j V_m {\boldsymbol y} _j  \| }{\| \sum_{j=1}^{n_0} w_j V_m {\boldsymbol y} _j  \|} 
    = \frac{ \| V_m  \sum_{j=1}^{2n_0} w_j {\boldsymbol y}_j  \| }{\| V_m \sum_{j=1}^{n_0} w_j  {\boldsymbol y} _j  \|}.
\end{equation} 
Note that 
\begin{equation}
    \left \| V_m \sum_{j=1}^{n_0} w_j  {\boldsymbol y}_j  \right \|^2 
    = \left(\sum_{j=1}^{n_0} w_j  {\boldsymbol y}_j\right)^{T} V_m^{T} V_m \sum_{j=1}^{n_0} w_j  {\boldsymbol y} _j 
    = \left \| \sum_{j=1}^{n_0} w_j  {\boldsymbol y}_j \right \|^2
\end{equation}
since $V_m^{T}V_m $ is the identity matrix. Dropping $V_m$ in \eqref{eq:reduced}, we define a new indicator
\begin{equation}\label{DeltaSm}
\tilde{I}_R = \frac{\|\sum_{j=1}^{2n_0} w_j {\boldsymbol y}_j \|}{\|\sum_{j=1}^{n_0} w_j  {\boldsymbol y}_j \|}.
\end{equation}
As a consequence, there is no need to store $V_m$'s ($n \times m$ matrices) but to store much smaller $m \times m$ ($m = O(1)$) matrices $H_m$'s.

As before, we use a threshold to decide whether or not eigenvalues exist in $R$. 
From \eqref{IndicatorPf}, if there are no eigenvalues in $R$, the indicator $I_R = O(e^{-C n_0})$.
In the experiments, we take $n_0=4$. Assume that $C=1$, we would have that $I_R \approx 0.018$. It is reasonable to
take $\delta_0 = 1/20$ as the threshold. The choice is ad-hoc. 
Nonetheless, the numerical examples show that the choice is robust. 

\begin{definition}
A (square) region $R$ is admissible if $I_R > \delta_0$.
\end{definition}
\begin{remark}
In practice, a region which is smaller than $h_0$ and not resolvable with respect to $(\sigma, \epsilon_0)$ is taken to be admissible.
\end{remark}

\subsection{Speedup the Computation of Indicators}
To check if a linear system \eqref{linearsys} can be solved effectively using a Krylov space $K^\sigma_{m}(M; {\boldsymbol b})$,
one need to compute the residual \eqref{eq:error} for many $z_j$'s. In the following, we propose a fast method.
First rewrite \eqref{eq:many} as
\begin{equation}\label{zjHm}
\left(\frac{1}{\sigma -z_j}I+H_m \right) {\boldsymbol y}_j = \frac{\beta}{\sigma -z_j} {\boldsymbol e}_1.
\end{equation}
Assume that $H_m$ has the following eigen-decomposition $H_m = PDP^{-1}$ where 
\[
D=\text{diag}\{ \lambda_1, \lambda_1, \ldots, \lambda_m \}.
\]
Then \eqref{zjHm} can be written as
\[
P\left(\frac{1}{\sigma -z_j}I+D \right)P^{-1}  {\boldsymbol y}_j = \frac{\beta}{\sigma -z_j} {\boldsymbol e}_1,
\]
whose solution is simply
\begin{eqnarray*}
{\boldsymbol y}_j  &=& P \left(\frac{1}{\sigma-z_j}I+D\right)^{-1}P^{-1} \frac{\beta}{\sigma- z_j}{\boldsymbol e}_1 \\
			&=& P\left(I + (\sigma-z_j)D\right)^{-1} P^{-1} {\boldsymbol e}_1.
\end{eqnarray*}
Hence
\begin{eqnarray}\nonumber
{\boldsymbol e}_m^{T} {\boldsymbol y}_j &=& {\boldsymbol e}_m^{T}  P\left(I + (\sigma-z_j)D\right)^{-1} P^{-1} {\boldsymbol e}_1 \\
\label{rmLc1} &=& {\boldsymbol r}_m \Lambda {\boldsymbol c}_1,
\end{eqnarray}
where ${\boldsymbol r}_m$ is the last row of $P$, ${\boldsymbol c}_1$ is the first column of $P^{-1}$, and
\[
 \Lambda = \text{diag}\left\{\frac{1}{1+(\sigma-z_j) \lambda_1}, \frac{1}{1+(\sigma-z_j) \lambda_2},\ldots, \frac{1}{1+(\sigma-z_j) \lambda_m}\right\}.
\]
In fact, this further reduces the memory requirement since only three $m\times 1$ vectors, ${\boldsymbol r}_m$, ${\boldsymbol c}_1$,
and $ \Lambda$ are stored for each shift $\sigma$.

\subsection{Multilevel Technique}\label{MT}
Now we propose a multilevel technique, which is more efficient and suitable for parallelization.
In SIM-M, the following strategy is employed.

At level $1$, $R$ is divided uniformly into smaller squares $R^1_j, j=1, \ldots, N^1$. Collect all quadrature points $z^1_j$'s and
solve the linear systems \eqref{linearsys} accordingly. The indicators of $R^1_j$'s are computed and squares containing eigenvalues
are chosen. Indicators of the resolvable squares are computed. Squares containing eigenvalues are subdivided into smaller square.
Squares that are not resolvable are also subdivided into smaller squares. These squares are left to the next level.
At level $2$, the same operation is carried out. The process stops at level $K$ when.
the size of the squares is smaller than the given precision $h_0$. 

\subsection{Multiplicities of Eigenvalues}\label{ME}
The first two members of SIMs only output the eigenvalues. A function to find the multiplicities of the eigenvalues can be integrated into SIM-M. 

\begin{definition}
An eigenvalue $\lambda$ is said to be resolved by a shift $\sigma$ if the small square at level $K$ containing $\lambda$ is
resolvable using the Krylov subspace $K_m^\sigma$.
\end{definition}

When the eigenvalues are computed, a mapping from the set of eigenvalues $\Lambda$ to the set of shifts $\Sigma$ is also established. 
Hence, for a shift $\sigma$, one can find the set of all eigenvalues that are resolved by $\sigma$, denoted by
\[
\Lambda_\sigma = \{ \lambda_1, \ldots, \lambda_n\}.
\]
For $k$ random vectors ${\boldsymbol f}_1, \ldots, {\boldsymbol f}_k$, generate $k$ Krylov subspaces $K_m^\sigma(M, {\boldsymbol b}_i), i=1, \ldots, k$.
For each $\lambda \in \Lambda_\sigma $, compute the spectral projections of ${\boldsymbol f}_1, \ldots, {\boldsymbol f}_k$
using the above Krylov subspaces.
Then the number of significant singular values of the matrix $[P{\boldsymbol f}_1, \ldots, P{\boldsymbol f}_k]$ is the multiplicity of $\lambda$.

\begin{remark}
In fact, the associated eigenvectors can be obtained with little extra cost by adding more quadrature points.
However, it can be expected that it needs a lot of more time and memory to find the multiplicities since more Krylov subspaces are generated.
\end{remark}

\subsection{Algorithm for SIM-M}
Now we are ready to present the new algorithm SIM-M.
\begin{itemize}
\item[] \text{SIM-M}$(A, B, R, {\boldsymbol f}, h_0, \epsilon, \delta_0, m, n_0)$
\item[] {\bf Input:} 
	\begin{itemize}
		\item $A, B$: $n \times n$ matrices 
		\item $R$: search region in $\mathbb C$
		\item ${\boldsymbol f}$: a random vector
		\item $h_0$: precision 
		\item $\epsilon$: residual tolerance
		\item $\delta_0$: indicator threshold
		\item $m$: size of Krylov subspace
		\item $n_0$: number of quadrature points
	\end{itemize}
\item[] {\bf Output:} 
	\begin{itemize}
	\item generalized eigenvalues $\lambda$'s inside $R$
	\end{itemize}
\item[1.] use the center of $R$ as the first shift and generate the associated Krylov subspaces.
\item[2.] pre-divide $R$ into small squares of size $h_0$: $R_j, j=1, \ldots, J$ (these are selected squares at the initial level).
\item[3.] for $j=1:J$ do
	\begin{itemize}
	\item For all quadrature points for $R_j$, check if the related linear systems can be solved using any one of the existing Krylov subspaces
		up to the given residual $\epsilon_0$. If yes, associate $R_j$ with that Krylov subspace. Otherwise, set the shift to be the center of
		$R_j$ and construct a Krylov subspace.
	\end{itemize}
\item[4.] calculate the number of the levels, denoted by $K$, needed to reach the precision $h_0$.
\item[5.] for $k=1:K$
	\begin{itemize}
		\item for each selected square $R_j^k$ at level $k$, check if $R_j^k$ is resolvable.
			\begin{itemize}
			\item if $R_j^k$ is resolvable, compute the indicator for $R_j^k$ and mark it when the indicator is larger than $\delta_0$, 
				i.e., $R_j^k$ contains eigenvalues.
			\item if $R_j^k$ is not solvable, mark $R_j^k$ and leave it to next level.
			\end{itemize}
		\item divide marked squares into four squares uniformly and move to next level.
	\end{itemize}
\item[6.] post-processing the marked squares at level $K$, merge eigenvalues when necessary, show warnings if there exist unsolvable squares.
\item[7.] output eigenvalues.
\end{itemize}

In the implementation, we choose $m=50$. Similar values such as $m=30$ do not change the performance significantly.
The indicator threshold is set to be $\delta_0 = 1/20$ as discussed in Section 3.1. The number of quadrature points is $n_0 = 8$, which is effective for the examples. 
The choices of these parameters affect the efficiency and robustness of the algorithm in a subtle way and
deserve more study for different problems.

\section{Numerical Examples}\label{NE}
We show some examples for SIM-M. All the test matrices are from the University of Florida Sparse Matrix Collection \cite{DavisHu2011ACMTOMS}
except the last example. The computations are done using MATLAB R2017a on a MacBook Pro with 16 GB memory and a 3-GHz Intel Core i7 CPU.

\subsection{Directed Weighted Graphs}
The first group contains four non-symmetric matrices, HB/gre\_115, HB/gre\_343, HB/gre\_512, HB/gre\_1107. 
These matrices represent directed weighted graphs. 

\begin{figure}[h!]
\label{grePlots}
\begin{minipage}{0.50\linewidth}
\includegraphics[angle=0, width=\textwidth]{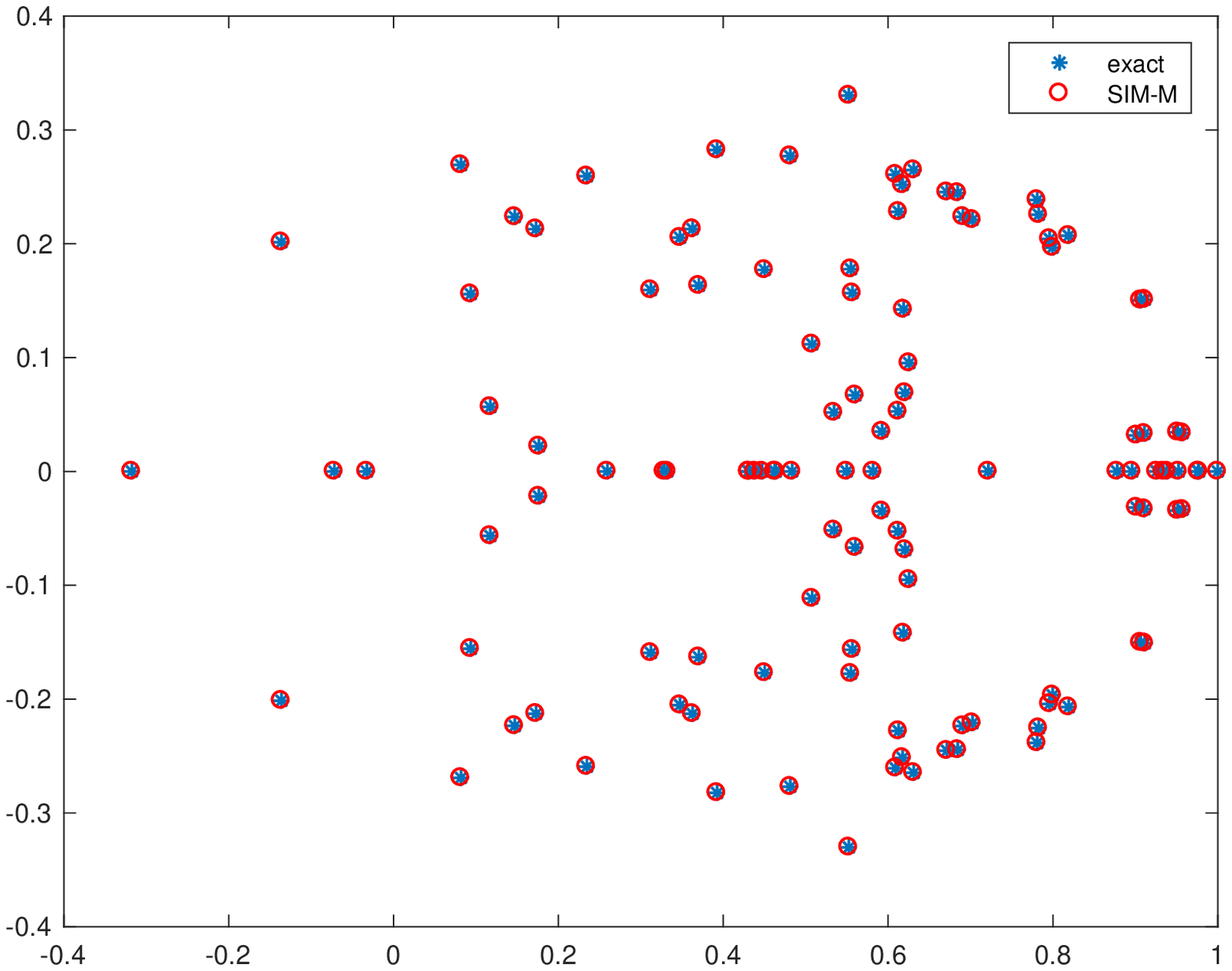}
\vspace{-0.5cm}
\begin{center}
(a)
\end{center}
\end{minipage}
\begin{minipage}{0.50\linewidth}
\includegraphics[angle=0, width=\textwidth]{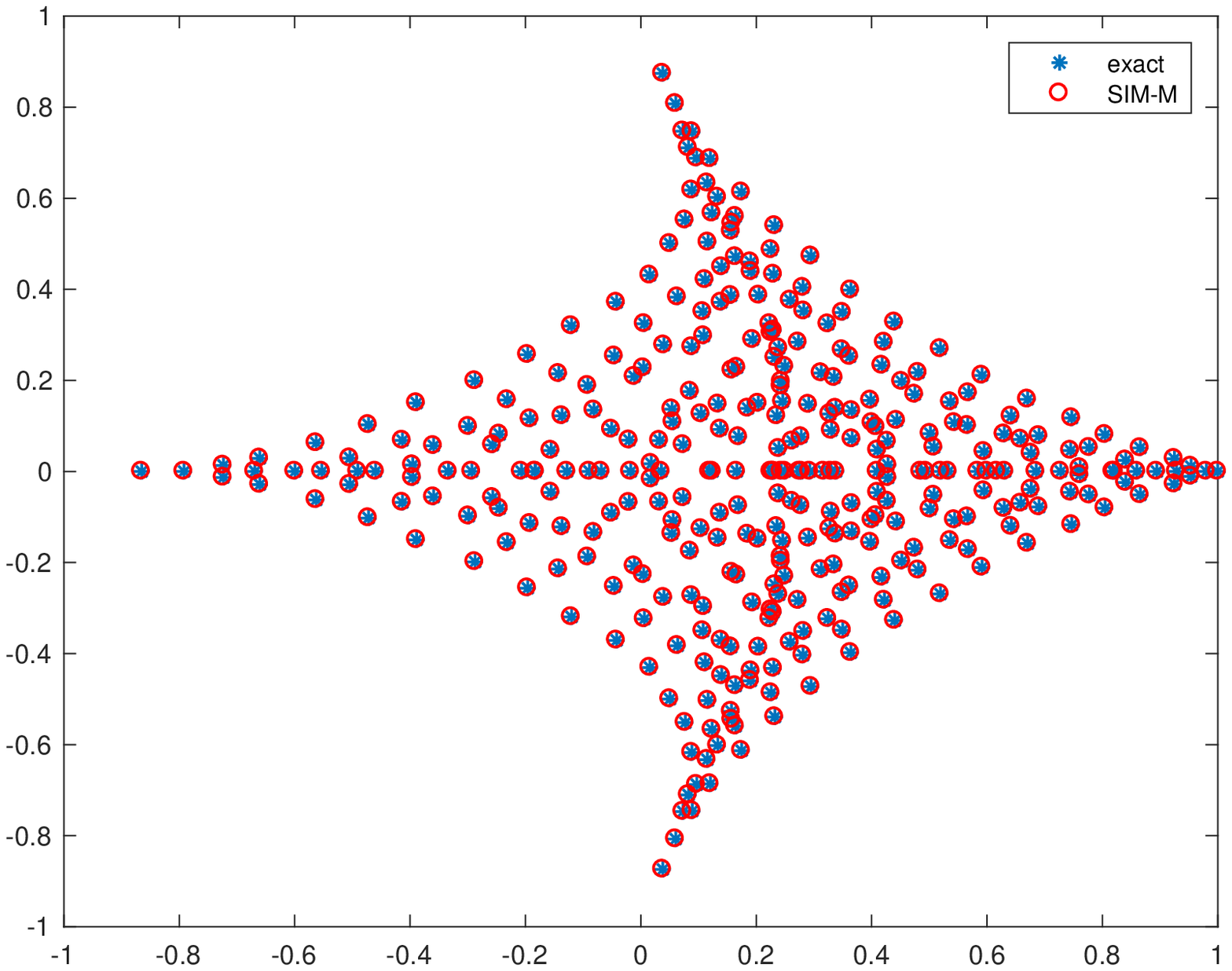}
\vspace{-1cm}
\begin{center}
(b)
\end{center}
\end{minipage}
\begin{minipage}{0.50\linewidth}
\includegraphics[angle=0, width=\textwidth]{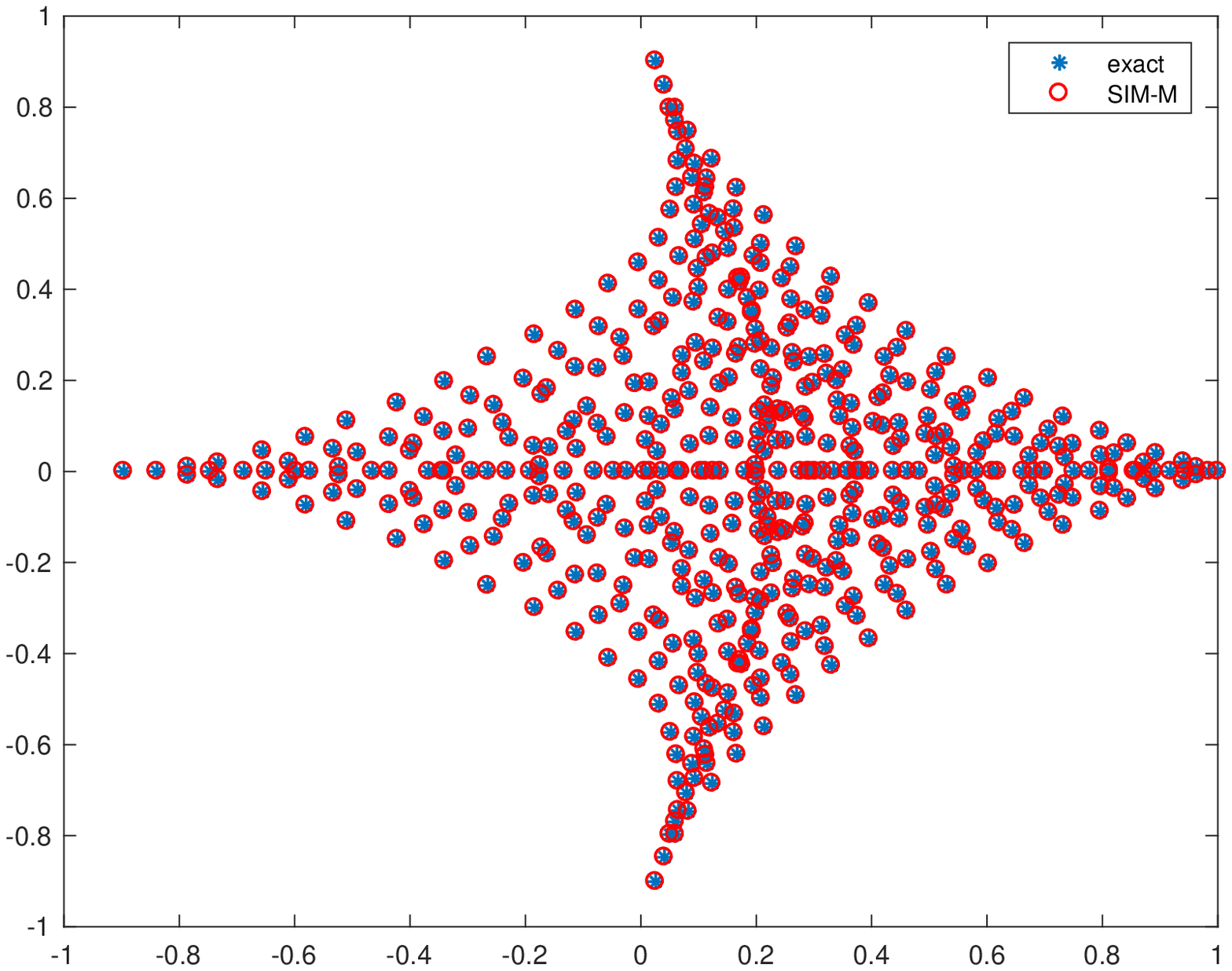}
\vspace{-1cm}
\begin{center}
(c)
\end{center}
\end{minipage}
\begin{minipage}{0.50\linewidth}
\includegraphics[angle=0, width=\textwidth]{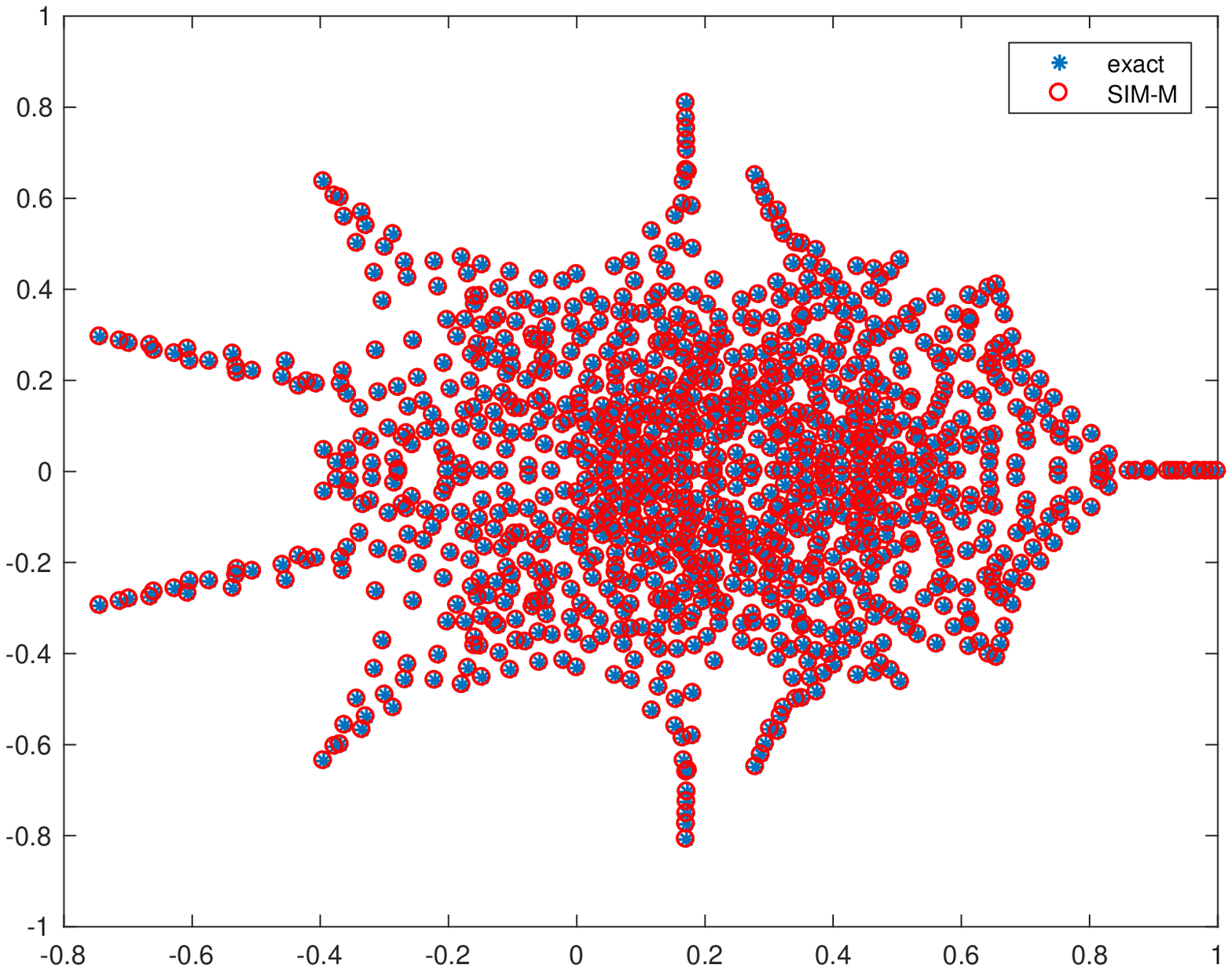}
\vspace{-1cm}
\begin{center}
(d)
\end{center}
\end{minipage}
\caption{Eigenvalues computed by SIM-M and Matlab {\it eig} coincide. (a) HB/gre\_115. (b) HB/gre\_343. 
(c): HB/gre\_512. (d): HB/gre\_1107.}
\end{figure}

\begin{table}[h!]
\caption{Time (in second) used for all eigenvalues by SIM-M.}
\label{gre}
\centering
\begin{tabular}{l|r|r|r|r}
\hline
N (size of the matrix) & 115& 343&  512& 1107\\ \hline
T (time in seconds) & 3.4141s & 10.2917s & 14.7461s& 40.2252s\\ \hline
T/N &0.0297  &  0.0300  &  0.0288  &  0.0363 \\ \hline
\hline
\end{tabular}
\end{table}
We compute all eigenvalues using SIM-M in Table~\ref{gre}. The first row represents sizes of the four matrices. The second row
shows the CPU times (in seconds) used by SIM-M. The numbers in the third row are the ratios of the seconds used by SIM-M and the sizes of the matrices,
i.e., the average time to compute one eigenvalue.
It seems that the ratio is stable for matrices of different sizes.
In Fig.~\ref{grePlots}, we show the eigenvalues computed by SIM-M and by Matlab {\it eig}, which coincide each other.

\subsection{A Quantum Chemistry Problem}
\begin{figure}[h!]
\label{QC2534}
\begin{minipage}{0.50\linewidth}
\includegraphics[angle=0, width=\textwidth]{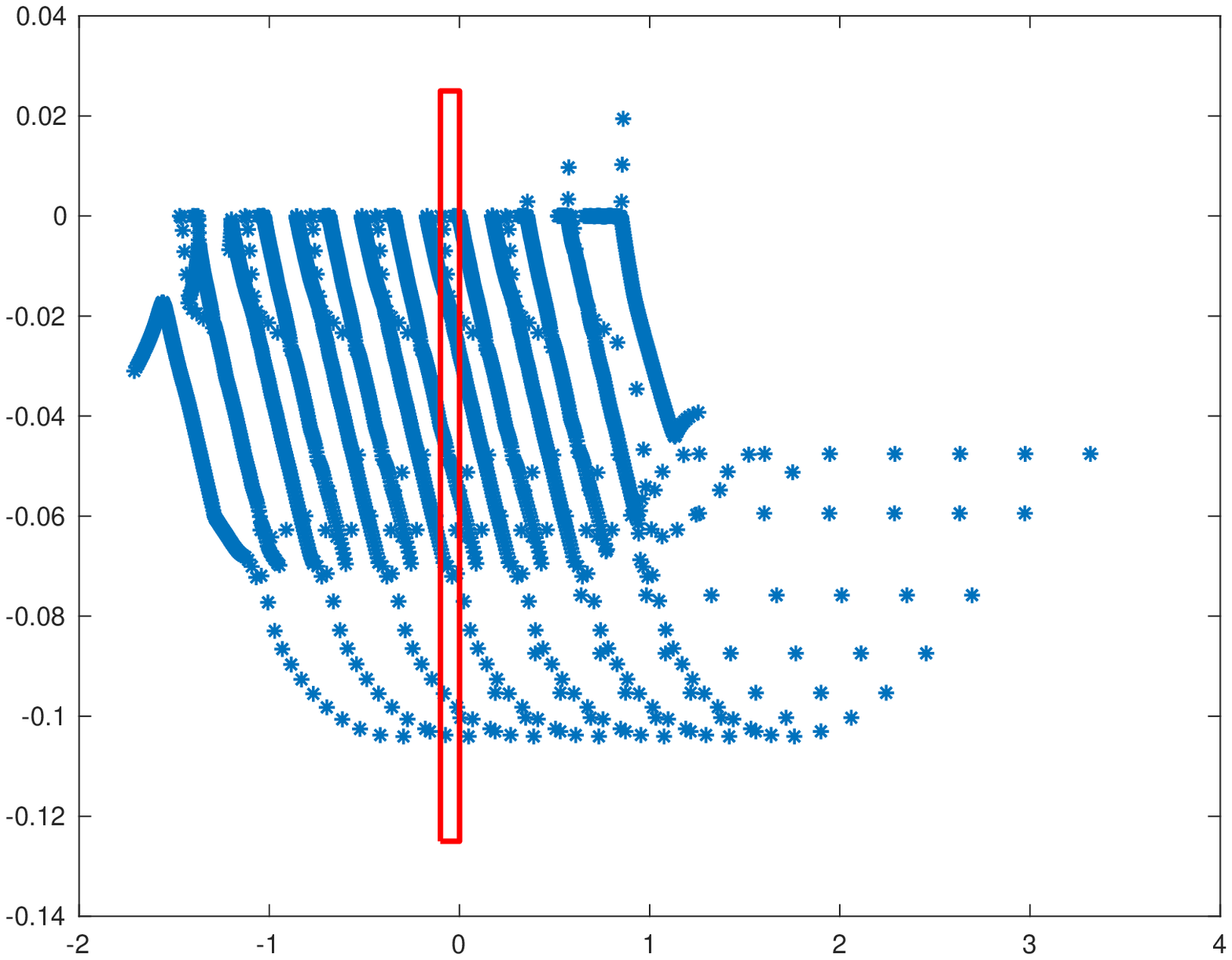}
\vspace{-0.5cm}
\begin{center}
(a)
\end{center}
\end{minipage}
\begin{minipage}{0.50\linewidth}
\includegraphics[angle=0, width=\textwidth]{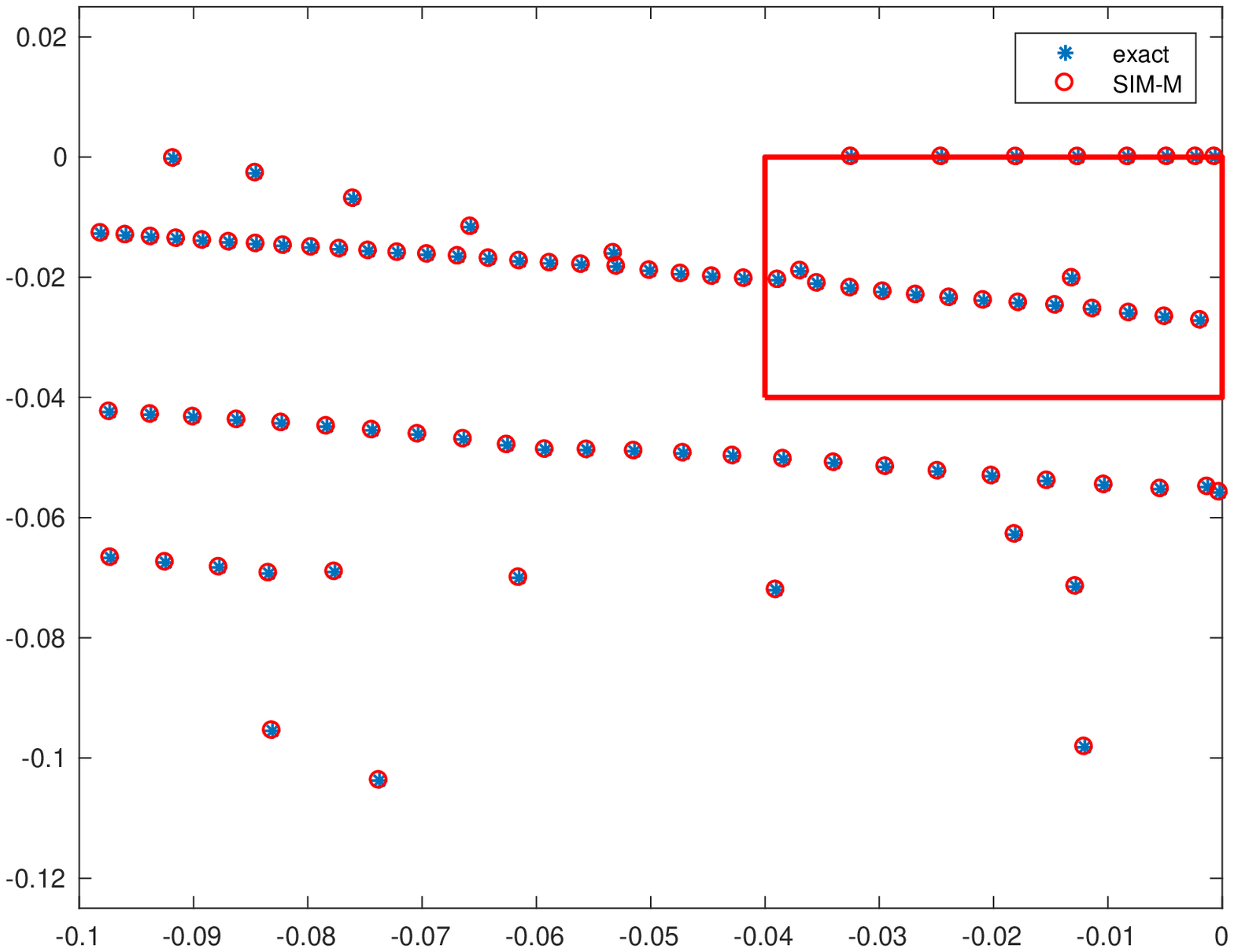}
\vspace{-1cm}
\begin{center}
(b)
\end{center}
\end{minipage}
\begin{minipage}{0.50\linewidth}
\includegraphics[angle=0, width=\textwidth]{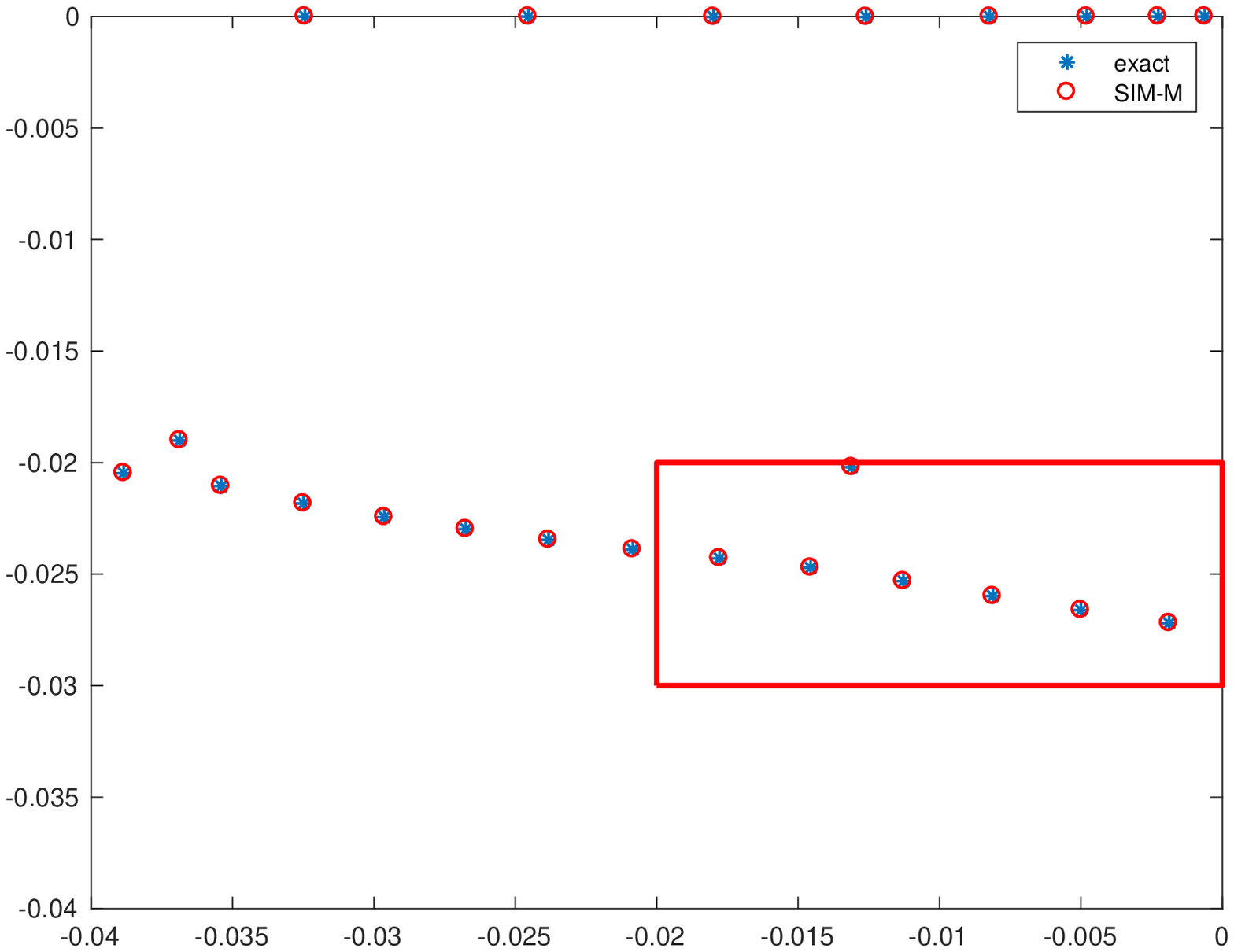}
\vspace{-1cm}
\begin{center}
(c)
\end{center}
\end{minipage}
\begin{minipage}{0.50\linewidth}
\includegraphics[angle=0, width=\textwidth]{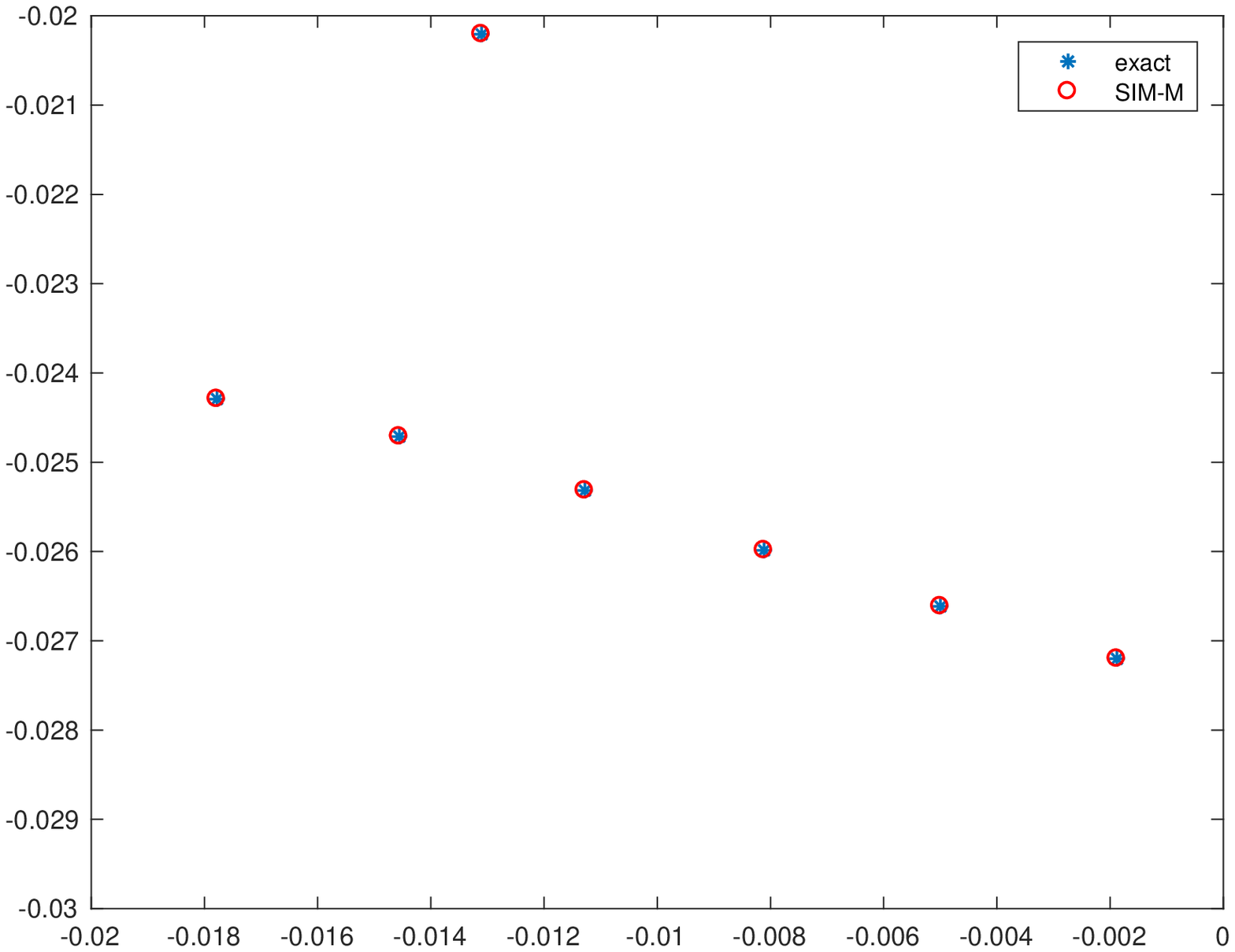}
\vspace{-1cm}
\begin{center}
(d)
\end{center}
\end{minipage}
\caption{QC2534. (a): Full spectrum by Matlab {\it eig} (the rectangle is $R_1$). 
(b): Eigenvalues by SIM-M in $R_1$ (the rectangle is $R_2$).
(c): Eigenvalues by SIM-M in $R_2$ (the rectangle is $R_3$).
(d): Eigenvalues by SIM-M in $R_3$.
}
\end{figure}

The second example, Bai/qc2534, is a sparse $2534 \times 2534$ matrix from modeling H2+ in an electromagnetic field.
The full spectrum, computed by Matlab {\it eig}, is shown in Fig.~\ref{QC2534}(a), in which
the red rectangle is $R_1=[-0.1, 0] \times [-0.125, 0.025]$. In Fig.~2(b), the eigenvalues are computed by SIM-M in $R_1$, which
coincide with those computed by Matlab {\it eig}. The red rectangle is Fig.~2(b) is $R_2=[-0.04, 0]\times [-0.04, 0]$.
Eigenvalues in $R_2$ computed by SIM-M are shown in Fig.~2(c). The rectangle in Fig.~2(c) is $R_3=[-0.02, 0]\times [-0.03, -0.02]$.
Eigenvalues in $R_3$ computed by SIM-M are shown in Fig.~2(d).

The second row of Table~\ref{T2534} shows that there are $88$, $23$ and $7$ eigenvalues in $R_1, R_2$ and $R_3$, respectively.
The third row shows the time used by SIM-M to compute all eigenvalues in $R_1, R_2$ and  $R_3$. 
The fourth row shows the average time to compute one eigenvalue, which seems to be consistent.
\begin{table}[h!]
\caption{Time (in second) used by SIM-M for different regions.}
\label{T2534}
\centering
\begin{tabular}{l|r|r|r}
\hline
& $R_1$&  $R_2$& $R_3$\\ \hline
N (\# of eigenvalues) & 88 & 23& 7\\  \hline
T (time in seconds) & 14.7445s & 3.7005s& 0.54645s\\  \hline
T/N & 0.1676  &  0.1609 & 0.0781\\  
\hline
\end{tabular}
\end{table}

\subsection{DNA Electrophoresis}
The third example is a $39,082 \times 39,082$ matrix, vanHeukelum/cage11, arising from DNA electrophoresis.
We consider a series of nested domains 
\begin{eqnarray*}
&&R_1=[0.230, 0.270] \times [-0.0005, 0.0005],\\
&&R_2=[0.250, 0.270] \times [-0.0005, 0.0005],\\
&&R_3=[0.250, 0.260] \times [-0.0005, 0.0005],\\
&&R_4=[0.254, 0.256] \times [-0.0005, 0.0005].
\end{eqnarray*}

In Table~\ref{cage11}, the time and number of eigenvalues found in each domain are shown.
Again, the average time to compute one eigenvalue is stable.
\begin{table}[h!]
\caption{Time (in second) used by SIM-M for different regions.}
\label{cage11}
\centering
\begin{tabular}{l|c|c|c|c}
\hline
& $R_1$&  $R_2$& $R_3 $& $R_4$\\ \hline\hline
N (\# of eigenvalues) &$105$&  $31$& $31$& $8$\\ \hline
T (time in seconds) & 588.3552s& 299.4242s & 214.0637s& 47.8098s\\ \hline
T/N & 5.6034  &  9.6588  &  6.9053  &  5.9762\\ 
\hline
\end{tabular}
\end{table}
\begin{remark}
Note that it is not possible to use Matlab {\it eig} to find all eigenvalues
due the memory constraint. However, SIM-M does not have this limitation.
In fact, numerical results in the above two subsections indicate that a parallel version of SIM-M 
has the potential to be faster than the classical methods. 
\end{remark}

\subsection{Quantum States in Disordered Media}
The test matrices are sparse and symmetric arising from localized quantum states in random or disordered media \cite{Arnold2016}.
We would like to use this example to show that the method can treat rather large problems on a laptop.
The matrices $A$ and $B$ are of $1,966,080 \times 1,966,080$. We consider three nested domains given by
\begin{eqnarray*}
&& R_1=[0.00, 0.60] \times [-0.05, 0.05],\\
&& R_2=[0.00, 0.50] \times [-0.05, 0.05], \\
&& R_3=[0.00, 0.40] \times [-0.05, 0.05], 
\end{eqnarray*}

In Table~\ref{AB65536}, time and number of eigenvalues in each domain are shown. 
Again, we observe that the average time to compute one eigenvalue is stable.
\begin{table}[h!]
\caption{Time (in second) used by SIM-M for different regions.}
\label{AB65536}
\centering
\begin{tabular}{l|c|c|c}
\hline
& $R_1$&  $R_2$& $R_3$\\ \hline \hline
N (\# of eigenvalues)  & $36$& $7$& $3$\\ \hline
T (time in seconds)   & 573.1088s&112.1876s & 58.9957s\\ \hline
T/N & 15.9197  & 16.0268   &19.6652 \\
\hline
\end{tabular}
\end{table}

\section{Conclusions and Future Work}\label{Con}
Given a region on the complex plane, SIMs first compute an indicator, which is used to test if the region contains eigenvalues.
Then the region is subdivided and tested until all the eigenvalues are isolated with a specified precision. Hence SIMs can be viewed as
a bisection technique. 

We propose an improved version SIM-M to compute many eigenvalues of large matrices. Several examples are presented for demonstrations.
However, to make the method practically competitive, a parallel implementation on super computers is necessary.
Currently, SIMs use the spectral projection to compute the indicators. Other ways to define the indicators should be investigated in future.

\end{document}